\documentclass[11pt]{amsart}
\usepackage{amsmath,amsthm,amsfonts,amssymb,times}
\usepackage{latexsym}

\usepackage{versions}
\usepackage{mathrsfs} 

\usepackage[notcite,notref]{showkeys} 

\numberwithin{equation}{section}  

\DeclareMathAlphabet{\curly}{U}{rsfs}{m}{n}  

\theoremstyle{remark}

\theoremstyle{plain}

\newtheorem{lem}{Lemma}[section]
\newtheorem{thm}{Theorem}
\newtheorem{cor}{Corollary}

\newcommand{\Z}{\mathbb{Z}}
\newcommand{\R}{\mathbb{R}}
\newcommand{\E}{\mathbb{E}}   
\newcommand{\PR}{\mathbb{P}}  

\makeatletter
\renewcommand{\pmod}[1]{\allowbreak\mkern7mu({\operator@font mod}\,\,#1)}
\makeatother

\newcommand{\bal}{\[\begin{aligned}}
\newcommand{\eal}{\end{aligned}\]}

\newcommand{\be}{\begin{equation}}
\newcommand{\ee}{\end{equation}}

\renewcommand{\a}{\ensuremath{\alpha}}

\newcommand{\del}{\ensuremath{\delta}}
\newcommand{\eps}{\ensuremath{\varepsilon}}

\renewcommand{\le}{\leqslant}
\renewcommand{\leq}{\leqslant}
\renewcommand{\ge}{\geqslant}
\renewcommand{\geq}{\geqslant}

\newcommand{\order}{\asymp}      
\renewcommand{\(}{\left(}
\renewcommand{\)}{\right)}

\newcommand{\pfrac}[2]{\left(\frac{#1}{#2}\right)}  

\newcommand{\ba}{\ensuremath{\mathbf{a}}}

\newcommand{\asym}{\sim}   

\newcommand{\PP}{\mathcal{P}}
\newcommand{\QQ}{\mathcal{Q}}
\newcommand{\cS}{\mathcal{S}}
\newcommand{\cR}{\mathcal{R}}

\renewcommand{\mod}{\bmod}  

\begin{document}

\title{Chains of large gaps between primes}

\author{Kevin Ford}
\address{Department of Mathematics\\ 1409 West Green Street \\ University
of Illinois at Urbana-Champaign\\ Urbana, IL 61801\\ USA}
\email{ford@math.uiuc.edu}



\author{James Maynard}
\address{Mathematical Institute\\
Radcliffe Observatory Quarter\\
Woodstock Road\\
Oxford OX2 6GG\\
England }
\email{james.alexander.maynard@gmail.com}

\author{Terence Tao}
\address{Department of Mathematics, UCLA\\
405 Hilgard Ave\\
Los Angeles CA 90095\\
USA}
\email{tao@math.ucla.edu}

\begin{abstract} Let $p_n$ denote the $n$-th prime, and for any $k \geq 1$ and sufficiently large $X$, define the quantity
$$ G_k(X) := \max_{p_{n+k} \leq X} \min( p_{n+1}-p_n, \dots, p_{n+k}-p_{n+k-1} ),$$
which measures the occurrence of chains of $k$ consecutive large gaps of primes.  Recently, with Green and Konyagin, the authors showed that
\[ G_1(X) \gg \frac{\log X \log \log X\log\log\log\log X}{\log \log \log X}\]
for sufficiently large $X$.  In this note, we combine the arguments in that paper with the Maier matrix method to show that
\[ G_k(X) \gg \frac{1}{k^2} \frac{\log X \log \log X\log\log\log\log X}{\log \log \log X}\]
for any fixed $k$ and sufficiently large $X$.  The implied constant is effective and independent of $k$.
\end{abstract}

\maketitle

\setcounter{tocdepth}{1}
\tableofcontents

\section{Introduction}

Let $p_n$ denote the $n^{\operatorname{th}}$ prime, and for any $k \geq 1$ and sufficiently large $X$, let
$$ G_k(X) := \max_{p_{n+k} \leq X} \min( p_{n+1}-p_n, \dots, p_{n+k}-p_{n+k-1} ),$$
denote the maximum gap between $k$ consecutive primes less than $X$.  The quantity $G_1(X)$ has been extensively studied.  The prime number theorem implies that
\[ G_1(X) \geq (1 + o(1)) \log X,\]
with the bound being successively improved in many papers \cite{back}, \cite{brauer}, \cite{West}, \cite{erdos-gaps}, \cite{R1}, \cite{schonhage}, \cite{rankin-1963},\cite{MP}, \cite{P}, \cite{maynard-large}, \cite{FGKT}, \cite{FGKMT}.  The best lower bound currently is\footnote{As usual in the subject, $\log_2 x := \log \log x$, $\log_3 x := \log \log \log x$, and so on.  The conventions for asymptotic notation such as $\ll$ and $o()$ will be defined in Section \ref{not-sec}.}
\[
G_1(X) \gg \frac{\log X \log_2 X \log_4 X}{\log_3 X},
\]
for sufficiently large $X$ and an effective implied constant, due to \cite{FGKMT}.  This result may be compared against the conjecture $G_1(X) \asymp \log^2 X$ of Cram\'er \cite{Cra} (see also \cite{Gra}), or the upper bound $G_1(X) \ll X^{0.525}$ of Baker-Harman-Pintz \cite{BHP}, which can be improved to $G_1(X) \ll X^{1/2} \log X$ on the Riemann hypothesis \cite{Cra1920}.

Now we turn to $G_k(X)$ in the regime where $k \geq 1$ is fixed, and $X$ assumed sufficiently large depending on $k$.  Clearly $G_k(X) \leq G_1(X)$, and a naive extension of the probabilistic heuristics of Cram\'er \cite{Cra} suggest that $G_k(X) \asymp \frac{1}{k} \log^2 X$ as $X \to \infty$.  The first non-trivial bound on $G_k(X)$ for $k \geq 2$ was by Erd\H{o}s \cite{erdos-gaps}, who showed that 
$$ G_2(X) / \log X \to \infty$$
as $X \to \infty$.  Using what is now known as the Maier matrix method, together with the arguments of Rankin \cite{R1} on $G_1(X)$, Maier \cite{maier} showed that
$$ G_k(X) \gg_k \frac{\log X \log_2 X \log_4 X}{(\log_3 X)^2}$$
for any fixed $k \geq 1$ and a sequence of $X$ going to infinity.  Recently, by modifying Maier's arguments and using the more recent work on $G_1(X)$ in \cite{FGKT}, \cite{maynard-large}, this was improved by Pintz \cite{P-new} to show that 
$$ G_k(X) / \left( \frac{\log X \log_2 X \log_4 X}{(\log_3 X)^2} \right) \to \infty$$
for a sequence of $X$ going to infinity.

Our main result here is as follows.

\begin{thm}\label{chain}  Let $k \geq 1$ be fixed.  Then for sufficiently large $X$, we have
\[
G_k(X) \gg \frac{1}{k^2} \frac{\log X \log_2 X \log_4 X}{\log_3 X}.
\]
The implied constant is absolute and effective.
\end{thm}

Maier's original argument required one to avoid Siegel zeroes, which restricted his results to a sequence of $X$ going to infinity, rather than all sufficiently large $X$.  However, it is possible to modify his argument to remove the effect of any exceptional zeroes, which allows us to extend the result to all sufficiently large $X$ and also to make the implied constant effective.  The intuitive reason for the $\frac{1}{k^2}$ factor is that our method produces, roughly speaking, $k$ primes distributed ``randomly'' inside an interval of length about $\frac{\log X \log_2 X \log_4 X}{\log_3 X}$, and the narrowest gap between $k$ independently chosen numbers in an interval of length $L$ is typically of length about $\frac{1}{k^2} L$.

Our argument is based heavily on our previous paper \cite{FGKMT}, in particular using the hypergraph covering lemma from \cite[Corollary 3]{FGKMT} and the construction of sieve weights from \cite[Theorem 5]{FGKMT}.  The main difference is in refining the probabilistic analysis in \cite{FGKMT} to obtain good upper and lower bounds for certain sifted sets arising in the arguments in \cite{FGKMT}, whereas in the former paper only upper bounds were obtained.

We remark that in the recent paper \cite{BF}, the methods from \cite{FGKMT} were modified to obtain some information about the limit points of tuples of $k$ consecutive prime gaps normalized by factors slightly slower than $\frac{\log X \log_2 X \log_4 X}{\log_3 X}$; see Theorem 6.4 of that paper for a precise statement.

\subsection{Acknowledgments}

KF thanks the hospitality of the Institute of Mathematics and Informatics of the
Bulgarian Academy of Sciences.
The research of JM was conducted partly while  he was a CRM-ISM postdoctoral fellow at the Universit\'e de Montr\'eal, and partly while he was a Fellow by Examination at Magdalen College, Oxford.  

KF was supported by NSF grant DMS-1201442.
TT was supported by a Simons Investigator grant, the
James and Carol Collins Chair, the Mathematical Analysis \&
Application Research Fund Endowment, and by NSF grant DMS-1266164. 

The authors thank Tristan Freiberg for some corrections.

\subsection{Notational conventions}\label{not-sec}
  
In most of the paper, $x$ will denote an asymptotic parameter going to infinity, with many quantities allowed to depend on $x$.
The symbol $o(1)$ will stand for a quantity bounded in magnitude by $c(x)$, where $c(x)$ is a quantity that tends to zero as $x \to \infty$.
The same convention applies to the asymptotic
notation
$X \asym Y$, which means $X=(1+o(1))Y$, and $X \lesssim Y$, which means $X \leq (1+o(1)) Y$.  We use $X = O(Y)$, $X \ll Y$, and $Y \gg X$ to denote the claim that there is a constant $C>0$ such that $|X| \le CY$ throughout the domain of 
the quantity $X$.
We adopt the convention that $C$ is independent of any parameter
unless such dependence is indicated, e.g. by subscript such as  $\ll_k$. 
 In all of our estimates here, the constant $C$ will be effective (we will not rely on ineffective results such as Siegel's theorem).  If we can take the implied constant $C$ to equal $1$, we write $f = O_{\leq}(g)$ instead.  Thus for instance
$$ X = (1 + O_{\leq}(\eps)) Y$$
is synonymous with
$$ (1-\eps) Y \leq X \leq (1+\eps) Y.$$
Finally, we use $X \order Y$ synonymously with $X \ll Y \ll X$.

When summing or taking products over the symbol $p$, it is understood that $p$ is restricted to be prime.

Given a modulus $q$ and an integer $n$, we use $n \mod q$ to denote the congruence class of $n$ in $\Z/q\Z$.

Given a set $A$, we use $1_A$ to denote its indicator function, thus $1_A(x)$ is equal to $1$ when $x \in A$ and zero otherwise.  Similarly, if $E$ is an event or statement, we use $1_E$ to denote the indicator, equal to $1$ when $E$ is true and $0$ otherwise.  Thus for instance $1_A(x)$ is synonymous with $1_{x \in A}$.


We use $\# A$ to denote the cardinality of $A$, and
for any positive real $z$, we let $[z] := \{ n \in \mathbf{N}: 1 \leq
n \leq z \}$ denote the set of natural numbers up to $z$.

Our arguments will rely heavily on the probabilistic method.  Our random variables will mostly be discrete (in the sense that they take at most countably many values), although we will occasionally use some continuous random variables (e.g. independent real numbers sampled uniformly from the unit interval $[0,1]$).  As such, the usual measure-theoretic caveats such as ``absolutely integrable'', ``measurable'', or ``almost surely'' can be largely ignored by the reader in the discussion below.  We will use boldface symbols such as $\mathbf{X}$ or $\mathbf{a}$ to denote random variables (and non-boldface symbols such as $X$ or $a$ to denote deterministic counterparts of these variables).  Vector-valued random variables will be denoted in arrowed boldface, e.g. $\vec{\mathbf{a}} = (\mathbf{a}_p)_{p \in \PP}$ might denote a random tuple of random variables $\mathbf{a}_p$ indexed by some index set $\PP$.  

We write $\PR$ for probability, and $\E$ for expectation.   If $\mathbf{X}$ takes at most countably many values, we define the \emph{essential range} of $\mathbf{X}$ to be the set of all $X$ such that $\PR( \mathbf{X} = X )$ is non-zero, thus $\mathbf{X}$ almost surely takes values in its essential range.  We also employ the following conditional expectation notation.  If $E$ is an event of non-zero probability, we write
$$ \PR( F | E ) := \frac{\PR( F \wedge E )}{\PR(E)}$$
for any event $F$, and
$$ \E( \mathbf{X} | E ) := \frac{\E({\mathbf X} 1_E)}{\PR(E)}$$
for any (absolutely integrable) real-valued random variable ${\mathbf X}$.  If $\mathbf{Y}$ is another random variable taking at most countably many values, we define the conditional probability $\PR(F|\mathbf{Y})$ to be the random variable that equals $\PR(F|\mathbf{Y}=Y)$ on the event $\mathbf{Y}=Y$ for each $Y$ in the essential range of $\mathbf{Y}$, and similarly define the conditional expectation $\E( \mathbf{X} | \mathbf{Y} )$ to be the random variable that equals $\E( \mathbf{X} | \mathbf{Y} = Y)$ on the event $\mathbf{Y}=Y$.  We observe the idempotency property
\begin{equation}\label{idem}
\E ( \E({\mathbf X}|{\mathbf Y}) ) = \E \mathbf{X}
\end{equation}
whenever ${\mathbf X}$ is absolutely integrable and $\mathbf{Y}$ takes at most countably many values.

We will rely frequently on the following simple concentration of measure result.

\begin{lem}[Chebyshev inequality]\label{cheb}  Let ${\mathbf X}, {\mathbf Y}$ be independent random variables taking at most countably many values.  Let ${\mathbf Y}'$ be a conditionally independent copy of ${\mathbf Y}$ over ${\mathbf X}$; in other words, for every $X$ in the essential range of $\mathbf{X}$, the random variables ${\mathbf Y}, {\mathbf Y}'$ are independent and identically distributed after conditioning to the event $\mathbf{X}=X$.  Let $F( {\mathbf X}, {\mathbf Y})$ be a (absolutely integrable) random variable depending on ${\mathbf X}$ and ${\mathbf Y}$.  Suppose that one has the bounds
\begin{equation}\label{1-moment}
\E F( {\mathbf X}, {\mathbf Y} ) = \alpha + O(\eps \alpha)
\end{equation}
and
\begin{equation}\label{2-moment}
\E F( {\mathbf X}, {\mathbf Y} ) F( {\mathbf X}, {\mathbf Y}' )  = \alpha^2 + O(\eps \alpha^2)
\end{equation}
for some $\alpha, \eps > 0$ with $\eps = O(1)$.  Then for any $\theta > 0$, one has
\begin{equation}\label{conclusion}
\E( F( {\mathbf X}, {\mathbf Y} ) | {\mathbf X} ) = \alpha + O_{\leq}(\theta)
\end{equation}
with probability $1 - O( \frac{\eps\alpha^2}{\theta^2})$.
\end{lem}

\begin{proof} See \cite[Lemma 1.2]{FGKMT}.
\end{proof}

\section{Siegel zeroes}

As is common in analytic number theory, we will have to address the possibility of an exceptional Siegel zero.  As we want to keep all our estimates effective, we will not rely on Siegel's theorem or its consequences (such as the Bombieri-Vinogradov theorem).  Instead, we will rely on the Landau-Page theorem, which we now recall.  Throughout, $\chi$ denotes a Dirichlet character.

\begin{lem}[Landau-Page theorem]\label{page}  Let $Q \ge 100$.  Suppose that $L(s,\chi) = 0$ for some primitive character $\chi$ of modulus at most $Q$, and some $s = \sigma + it$.  Then either
$$ 1-\sigma \gg \frac{1}{\log(Q (1+|t|))},$$
or else $t=0$ and $\chi$ is a quadratic character $\chi_Q$, which is unique for any given $Q$.  Furthermore, if $\chi_Q$ exists, then its conductor $q_Q$ is square-free apart from a factor of at most $4$, and obeys the lower bound
$$ q_Q \gg \frac{\log^2 Q}{\log^2_2{Q}}.$$
\end{lem}

\begin{proof}  See e.g. \cite[Chapter 14]{Da}.  The final estimate follows from the classical bound $1-\beta\gg q^{-1/2}\log^{-2}{q}$ for a real zero $\beta$ of $L(s,\chi)$ with $\chi$ of modulus $q$. 
\end{proof}

We can then eliminate the exceptional character by deleting at most one prime
factor of $Q$.

\begin{cor}\label{page-cor}  Let $Q \ge 100$.  Then there exists a quantity $B_Q$ which is either equal to $1$ or is a prime of size
$$ B_Q \gg \log_2 Q$$
with the property that
$$ 1-\sigma \gg  \frac{1}{\log(Q (1+|t|))}$$
whenever $L(\sigma+it,\chi)=0$ and $\chi$ is a character of modulus at most $Q$ and coprime to $B_Q$.
\end{cor}

\begin{proof}  If the exceptional character $\chi_Q$ from Lemma \ref{page} does not exist, then take $B_Q := 1$; otherwise we take $B_Q$ to be the largest prime factor of $q_Q$.  As $q_Q$ is square-free apart from a factor of at most $4$, we have $\log q_Q \ll B_Q$ by the prime number theorem, and the claim follows.
\end{proof}

Next, we recall Gallagher's prime number theorem:

\begin{lem}[Gallagher's prime number theorem]\label{gallagher}  Let $q$ be a natural number, and suppose that $L(s,\chi) \neq 0$ for all characters $\chi$ of modulus $q$ and $s$ with $1-\sigma \leq \frac{\del}{\log(Q(1+it))}$, and some constant $\del>0$.  Then there is a constant $D \geq 1$ depending only on $\del$ such that
$$ \# \{ p \hbox{ prime}: p \leq x; p \equiv a\pmod{q} \} \gg \frac{x}{\phi(q) \log x}$$
for all $(a,q)=1$ and  $x \geq q^D$.
\end{lem}

\begin{proof} See \cite[Lemma 2]{maier}.
\end{proof}

This will combine well with Corollary \ref{page-cor} once we remove
 the moduli divisible by the (possible) exceptional prime $B_Q$.

\section{Sieving an interval}\label{sec:sieving}

We now give the key sieving result that will be used to prove Theorem \ref{chain}.

\begin{thm}[Sieving an interval]\label{sieve-thm} 
There is an absolute constants $c>0$ such that the following holds.
Fix $A \geq 1$ and $\eps > 0$, and let $x$ be sufficiently large depending on $A$ and $\eps$.
Suppose $y$ satisfies
\begin{equation}\label{ydef}
 y = c \frac{x \log x \log_3 x}{\log_2 x},
\end{equation}
and suppose that $B_0=1$ or that $B_0$ is a prime satisfying
$$ \log x \ll B_0 \leq x.$$
Then one can find a congruence class $a_p \mod p$ for each prime $p \leq x$, $p \neq B_0$ such that the sieved set 
$$ {\mathcal T} := \{ n \in [y] \backslash [x]: n \not\equiv a_p \pmod p \hbox{ for all } p \leq x, p \neq B_0 \}$$
obeys the following size estimates:
\begin{itemize}
\item (Upper Bound) One has
\begin{equation}\label{up-bound}
 \# {\mathcal T} \ll A \frac{x}{\log x}.
\end{equation}
\item (Lower Bound) One has
\begin{equation}\label{down-bound}
 \# {\mathcal T} \gg A \frac{x}{\log x}.
\end{equation}
\item (Upper bound in short intervals)  For any $0 \leq \alpha \leq \beta \leq 1$, one has
\begin{equation}\label{up-short}
 \# ({\mathcal T} \cap [\alpha y, \beta y]) \ll A (|\beta-\alpha|+\eps) \frac{x}{\log x}.
\end{equation}
\end{itemize}
\end{thm}

We remark that if one lowers $y$ to be of order $\frac{x \log x \log_3 x}{(\log_2 x)^2}$ rather than $\frac{x \log x \log_3 x}{\log_2 x}$, then this theorem is essentially \cite[Lemma 6]{maier}.  It is convenient to sieve $[y] \backslash [x]$ instead of $[y]$ for minor technical reasons (we will use the fact that the residue class $0 \mod p$ avoids all the primes in $[y] \backslash [x]$ whenever $p \leq x$).  The arguments in \cite{FGKMT} already can give much of this theorem, with the exception of the lower bound \eqref{down-bound}, which is the main additional technical result of this paper that is needed to extend the results of that paper to longer chains.

We will prove Theorem \ref{sieve-thm} in later sections.  In this section, we show how this theorem implies Theorem \ref{chain}.  Here we shall use the Maier matrix method, following the arguments in \cite{maier} closely (although we will use probabilistic notation rather than matrix notation).  Let $k \geq 1$ be a fixed integer, let $c_0>0$ be a small constant, and let $A \geq 1$ and $0 < \eps < 1/2$ be large and small quantities depending on $k$ to be chosen later.  

We now recall (a slight variant of) some lemmas from \cite{maier}.  

\begin{lem}\label{pxa}  There exists an absolute constant $D \geq 1$ such that, for all sufficiently large $x$, there exists a natural number $B_0$ which is either equal to $1$ or a prime, with
\begin{equation}\label{b0-est}
 \log x \ll B_0 \leq x,
\end{equation}
and is such that the following holds. If one sets $P :=P(x)/B_0$ (where we recall that $P(x)$ is the product of the primes up to $x$), then one has
\begin{equation}\label{zpx}
 \# \{ z \in [Z]: P z + a \hbox{ prime} \} \gg \frac{\log x}{\log Z} Z
\end{equation}
for all $Z \geq P^D$ and $a \in P$ coprime to $P$, and
\begin{equation}\label{zpy}
 \# \{ z \in [Z]: P z + a, P z + b \hbox{ both prime} \} \ll 
\pfrac{\log x}{\log Z}^2 Z
\end{equation}
for all $Z \geq P^D$ and all distinct $a,b \in [P]$ coprime to $P$.
\end{lem}

\begin{proof}   
We first prove \eqref{zpx}.   We apply Corollary \ref{page-cor} with $Q := P(x)$ to obtain a quantity $B_{P(x)}$ with the stated properties.  We set $B_0=1$ if $B_{P(x)} > x$, and $B_0 := B_{P(x)}$ otherwise. Then from Mertens' theorem we have \eqref{b0-est} if $B_0 \neq 1$.  From Corollary \ref{page-cor} and Lemma \ref{gallagher}, we then have
$$
 \# \{ z \in [Z]: P z + a \hbox{ prime} \} \gg \frac{PZ}{\phi(P) \log(PZ)} $$
for any $Z \geq P^D$ and a suitable absolute constant $D \geq 1$.  Note that $\log(PZ) \ll \log Z$.
From Mertens' theorem (and \eqref{b0-est}) we also have
\begin{equation}\label{ppx}
 \frac{P}{\phi(P)} \order \log x,
\end{equation}
and \eqref{zpx} follows.

Finally, the estimate \eqref{zpy} follows from standard upper bound sieves (cf. \cite[Lemma 3]{maier}).
\end{proof}

Now set $Z := P^D$ with $x$ and $D$ as in Lemma \ref{pxa}, and let $\mathbf{z}$ be chosen uniformly at random from $[Z]$.  Let $y$, ${\mathcal T}$ and $a_p \mod p$ be as in Theorem \ref{sieve-thm}.  By the Chinese remainder theorem, we may find $m \in [P]$ such that $m \equiv -a_p \pmod p$ for all $p \leq x$ with $p \neq B_0$.  Thus, $\mathbf{z}P + m+{\mathcal T}$ consists precisely of those elements of $\mathbf{z}P + m+[y] \backslash [x]$ that are coprime to $P$.  In particular, any primes that lie in the interval $\mathbf{z}P+m+[y] \backslash [x]$ lie in $\mathbf{z}P + m+{\mathcal T}$.

From \eqref{zpx} and Mertens' theorem we have
$$ \PR( \mathbf{z}P + m + a \hbox{ prime} ) \gg \frac{\log x}{x} $$
for all $a \in {\mathcal T}$ (we allow implied constants to depend on $D$).  Similarly, from \eqref{zpy} and Mertens' theorem we have
\begin{equation}\label{og}
 \PR( \mathbf{z}P + m + a, \mathbf{z}P(x) + m + b \hbox{ both
   prime} ) \ll \pfrac{\log x}{x}^2
\end{equation}
for any distinct $a,b \in {\mathcal T}$.
If we let $\mathbf{N}$ denote the number of primes in $\mathbf{z} P + m + {\mathcal T}$ (or equivalently, in $\mathbf{z} P + m + [y] \backslash [x]$), we thus have from \eqref{up-bound} and \eqref{down-bound} that
$$ \E \mathbf{N} \gg A $$
and
$$ \E \mathbf{N}^2 \ll A^2.$$
From this we see that with probability $\gg 1$, we have
\begin{equation}\label{ana}
 A \ll \mathbf{N} \ll A,
\end{equation}
where all implied constants are independent of $\eps$ and $A$.  (This is because the contribution to $\E \mathbf{N}$ when $\mathbf{N}$ is much larger than $A$ is much smaller than $A$.)

Next, if $0 \leq \alpha \leq \beta \leq 1$ and $\beta-\alpha \leq 2\eps$, then from \eqref{og}, \eqref{up-short} and the union bound we see that the probability that there are at least two primes in $\mathbf{z} P + m + [\alpha y, \beta y]$ is at most
$$ O\bigg( \(A \eps \frac{x}{\log x}\)^2 \pfrac{\log x}{x}^2 \bigg) = O( A^2 \eps^2 ).$$
Note that one can cover $[0,1]$ with $O(1/\eps)$ intervals of length
at most $2\eps$, with the property that any two elements $a,b$ of
$[0,1]$ with $|a-b| \leq \eps$ may be covered by at least one of these
intervals.  From this and the union bound, we see that the probability
that $\mathbf{z} P + m + [y] \backslash [x]$ contains two primes separated by at
most $\eps y$ is bounded by $O( \frac{1}{\eps} A^2 \eps^2 ) = O( A^2 \eps )$.  In particular, if we choose $\eps$ to be a sufficiently small multiple of $\frac{1}{A^2}$, we may find $z \in [Z]$ such that the interval $z P + m + [y] \backslash [x]$ contains $\gg A$ primes and has no prime gap less than $\eps y$.  If we choose $A$ to be a sufficiently large multiple of $k$, we conclude that
$$ G_k( Z P + m + y ) \geq \eps y \gg \frac{1}{k^2} y.$$
By Mertens' theorem, we have $Z P + m + y \ll \exp( O(x) )$, and Theorem \ref{chain} then follows from \eqref{ydef}.

It remains to prove Theorem \ref{sieve-thm}.  This is the objective of the remaining sections of the paper.

\section{Sieving a set of primes}\label{sec:initial}

Theorem \ref{sieve-thm} concerns the problem of deterministically
sieving an interval $[y] \backslash [x]$ of size \eqref{ydef} so that the sifted set
${\mathcal T}$ has certain size properties.  We use a variant of the
Erd\H{o}s-Rankin method to reduce this problem to a problem of
\emph{probabilistically} sieving a set $\QQ$ of \emph{primes} in $[y] \backslash [x]$,
rather than integers in $[y] \backslash [x]$. 

Given a real number $x \geq 1$, and a natural number $B_0$, define
\begin{equation}\label{zdef}
 z :=x^{\log_3 x/(4\log_2 x)},
\end{equation}
and introduce the three disjoint sets of primes
\begin{align}
\cS &:= \{ s\; \mbox{prime} : \log^{20} x < s \le z; s \neq B_0 \},\label{s-def}\\
\PP &:= \{ p \; \mbox{prime}: x/2 < p \leq x; p \neq B_0\},\label{p-def}\\
\QQ &:= \{ q \;\mbox{prime}: x < q \leq y; q \neq B_0\}\label{q-def}.
\end{align}
For residue classes $\vec a = (a_s \mod s)_{s\in S}$ and 
$\vec n = (n_p \mod p)_{p\in \PP}$, define the sifted sets
$$ S(\vec a) := \{ n \in \Z: n \not\equiv a_s \pmod s \hbox{ for all } s \in \cS \}$$
and likewise
$$ S(\vec n) := \{ n \in \Z: n \not\equiv n_p \pmod p \hbox{ for all } p \in \PP \}.$$

We reduce Theorem \ref{sieve-thm} to

\begin{thm}[Sieving primes]\label{sieve-primes}  Let $A \geq 1$ be a
  real number, let $x$ be sufficiently
  large depending on $A$, and suppose that $y$ obeys  \eqref{ydef}.
Let $B_0$ be a natural number.
Then there is a quantity 
	\begin{equation}\label{adef}
	A' \order A,
	\end{equation}
and some way to choose
the vectors $\vec{\mathbf{a}}=(\mathbf{a}_s \mod s)_{s \in \cS}$
and $\vec{\mathbf{n}}=(\mathbf{n}_p \mod p)_{p \in \PP}$  at random
(not necessarily independent of each other), such that
for any fixed $0 \leq \alpha < \beta \leq 1$ (independent of $x$), 
one has with probability $1-o(1)$ that
\begin{equation}\label{up-short-random}
 \#(\QQ \cap S( \vec{\mathbf a} ) \cap S( \vec{\mathbf n} ) \cap (\alpha y, \beta y]) \sim A' |\beta-\alpha| \frac{x}{\log x}.
\end{equation}
The $o(1)$ decay rates in the probability error and implied in the $\sim$ notation are allowed to depend on $A, \alpha,\beta$.
\end{thm}

In \cite[Theorem 2]{FGKMT}, a weaker version of this theorem was established in which $B_0$ was not present, and only the upper bound in \eqref{up-short-random} was proven.  Thus, the main new contribution of this paper is the lower bound in \eqref{up-short-random}.

We prove Theorem \ref{sieve-primes} in subsequent sections. In this section, we show how this theorem implies Theorem \ref{sieve-thm} (and hence Theorem \ref{chain}).  The arguments here are almost identical to those in \cite[\S 2]{FGKMT}.

Fix  $A\ge 1, 0<\eps\le 1$.  We partition $(0,1]$ into $O(1/\eps)$
  intervals $[\alpha_i,\beta_i]$ of length between $\eps/2$ and
  $\eps$.  Applying Theorem \ref{sieve-primes} with the pairs
  $(\a,\beta)=(\alpha_i,\beta_i)$ and the pair $(\a,\beta)=(0,1)$, and
  invoking a union bound (and the fact that $\eps$ is independent of $x$), we
see that if $x$ is sufficiently large (depending on $A,\eps$), there are $A', y$ obeying \eqref{adef}, \eqref{ydef} and tuples
of residue classes $\vec a = (a_s \mod s)_{s \in \cS}$ and 
$\vec n = (n_p \mod p)_{p  \in \PP}$ such that 
$$ \#(\QQ \cap S(\vec a) \cap S(\vec n)) \asym A' \frac{x}{\log x}$$
and
$$
 \#(\QQ \cap S(\vec a) \cap S(\vec n)) \cap (\alpha_i y, \beta_i y]) \ll A \eps \frac{x}{\log x}$$
for all $i$.  A covering argument then gives
$$
 \#(\QQ \cap S(\vec a) \cap S(\vec n) \cap [\alpha y, \beta y]) \ll A (|\beta-\alpha|+\eps) \frac{x}{\log x}$$
for any $0 \leq \alpha < \beta \leq 1$.
Now we extend the tuple $\vec a$ to a tuple $(a_p)_{p \leq x}$ of congruence classes $a_p \mod p$ for all primes $p \leq x$ by setting $a_p := n_p$ for $p \in \PP$ and $a_p := 0$ for $p \not \in \cS \cup \PP$, and consider the sifted set
$$ {\mathcal T} := \{ n \in [y] \backslash [x]: n \not\equiv a_p \pmod p \hbox{ for all } p \leq x \}.$$
The elements of ${\mathcal T}$, by construction, are not divisible by any prime in $(0,\log^{20} x]$ or in $(z,x/2]$, except possibly for $B_0$. 
 Thus, each element must either be a $z$-smooth number (i.e. a number with all prime factors at most $z$) times a power of $B_0$, or must consist of a prime greater than $x/2$, possibly multiplied by some additional primes that are all either at least $\log^{20} x$ or equal to $B_0$.  However, from \eqref{ydef} we know that $y=o(x\log x)$, and by hypothesis we know that $B_0 \gg \log x$.
Thus, we see that an element of ${\mathcal T}$ is either a $z$-smooth number times a power of $B_0$ or a prime in $\QQ$.  In the second case, the element lies in $\QQ  \cap S(\vec a) \cap S(\vec n)$.  Conversely, every element of $\QQ \cap S(\vec a) \cap S(\vec n)$ lies in ${\mathcal T}$.  Thus, ${\mathcal T}$ only differs from $\QQ  \cap S(\vec a) \cap S(\vec n)$ by a set $\cR$ consisting of $z$-smooth numbers in $[y]$ multiplied by powers of $B_0$. 

To estimate $\# \cR$, let
$$ u := \frac{\log y}{\log z},$$
so from \eqref{ydef}, \eqref{zdef} one has $u \asym 4 \frac{\log_2 x}{\log_3 x}$.  
The number of powers of $B_0$ in $[y]$ is $O(\log x)$.  By standard counts for smooth numbers (e.g. de Bruijn's theorem \cite{deB}) and \eqref{ydef}, we thus have
\begin{align*}
\# \cR  &\ll \log x \times y e^{-u\log u + O( u \log\log(u+2) ) } \\
&= \log x \times \frac{y}{\log^{4+o(1)} x} 
= o\( \frac{x}{\log x}\).
\end{align*}
Thus the contribution of $\cR$ to ${\mathcal T}$ is negligible for the purposes of establishing the bounds \eqref{up-bound}, \eqref{down-bound}, \eqref{up-short}, and Theorem \ref{sieve-thm} follows from \eqref{up-short-random}.

It remains to establish Theorem \ref{sieve-primes}.  This is the objective of the remaining sections of the paper.

\section{Using a hypergraph covering theorem}\label{sec:pip}

In the previous section we reduced matters to obtaining random residue classes $\vec{\mathbf{a}}$, $\vec{\mathbf{n}}$ such that the sifted set $\QQ \cap S(\vec{\mathbf a}) \cap S(\vec{\mathbf n})$ is small.  In this section we use a hypergraph covering theorem from \cite{FGKMT} to reduce the task to that of finding random residue classes $\vec{\mathbf n}$ that have large intersection with $\QQ \cap S(\vec{\mathbf a})$.  More precisely, we will use the following result:

\begin{thm}\label{packing-quant-cor}  Let $x\to\infty$.
Let $\PP'$, $\QQ'$ be sets of primes in $(x/2,x]$ and $(x, x\log x]$, respectively, with $\#\QQ' > (\log_2 x)^3$.
For each $p \in \PP'$, let $\mathbf{e}_p$ be a random subset of $\QQ'$
satisfying the size bound
\be\label{rbound}
\# \mathbf{e}_p \le r = O\( \frac{\log x \log_3 x}{\log_2^2 x} \) \qquad
(p\in \PP').
\ee
Assume the following:
\begin{itemize}
\item (Sparsity) For all $p \in \PP'$ and $q \in \QQ'$,
\begin{equation}\label{q-form-quant-cor}
\PR( q \in \mathbf{e}_p ) \le x^{-1/2 - 1/10}.
\end{equation}
\item (Uniform covering) For all but at most $\frac{1}{(\log_2 x)^2} \# \QQ'$ elements $q \in \QQ'$, we have
\begin{equation}\label{pje-size-bite-cor}
\sum_{p \in \PP'} \PR( q \in \mathbf{e}_p) = C + O_{\le}\pfrac{1}{(\log_2 x)^2}
 \end{equation}
for some quantity $C$, independent of $q$, satisfying
\begin{equation}\label{sigma}
\frac{5}{4} \log 5 \le C \ll 1.
\end{equation}
\end{itemize}
Then for any positive integer $m$ with
\begin{equation}\label{moo}
 m \le \frac{\log_3 x}{\log 5},
\end{equation}
we can find random sets $\mathbf{e}'_p\subseteq \QQ'$ for each $p \in \PP'$ such that
\[
\# \{ q \in \QQ':  q \not\in \mathbf{e}'_p \hbox{ for all } p \in \PP' \} \sim 5^{-m} \# \QQ'
\]
with probability $1-o(1)$.  More generally, for any $\QQ'' \subset \QQ'$ with cardinality at least $(\# \QQ')/\sqrt{\log_2 x}$, one has
\[
\# \{ q \in \QQ'':  q \not\in \mathbf{e}'_p \hbox{ for all } p \in \PP' \} \sim 5^{-m} \# \QQ''
\]
with probability $1-o(1)$.  The decay rates in the $o(1)$ and $\sim$ notation
are uniform in $\PP'$, $\QQ'$, $\QQ''$.
\end{thm}

\begin{proof} See \cite[Corollary 3]{FGKMT}.
\end{proof}

In view of the above result, we may now reduce Theorem \ref{sieve-primes} to the following claim.

\begin{thm}[Random construction]\label{sieve-primes-2}  Let $x$ be a sufficiently real number, let $B_0$ be a natural number and suppose $y$ satisfies \eqref{ydef}.  Then there is a quantity $C$ with
\begin{equation}\label{sigma-order}
C \order \frac{1}{c}
\end{equation}
with the implied constants independent of $c$, and some way to choose random vectors $\vec{\mathbf{a}}=(\mathbf{a}_s \mod s)_{s \in \cS}$
and $\vec{\mathbf{n}}=(\mathbf{n}_p)_{p \in \PP}$ of congruence classes $\mathbf{a}_s \mod s$ and integers $\mathbf{n}_p$, obeying the following axioms:
\begin{itemize}
\item For every $\vec a$ in the essential range of $\vec{\mathbf{a}}$, one has
$$ \PR( q \equiv \mathbf{n}_p \pmod p | \vec{\mathbf{a}} = \vec a) 
\le x^{-1/2 - 1/10}
$$
uniformly for all $p \in \PP$.
\item For fixed $0 \leq \alpha < \beta \leq 1$, we have with probability $1-o(1)$ that
\begin{equation}\label{treat}
 \# (\QQ \cap S(\vec{\mathbf{a}}) \cap [\alpha y, \beta y]) \sim 80 c |\beta-\alpha| \frac{x}{\log x} \log_2 x.
\end{equation}
\item Call an element $\vec a$ in the essential range of $\vec{\mathbf{a}}$
\emph{good} if, for all but at most 
$\frac{x}{\log x \log_2 x}$ elements $q \in \QQ \cap S(\vec{\mathbf{a}})$, one has
\be\label{good}
\sum_{p \in \PP} \PR( q \equiv \mathbf{n}_p \pmod p | \vec{\mathbf{a}}=\vec{a}) = C + O_{\le} \pfrac{1}{(\log_2 x)^2}.
\ee
Then  $\vec{\mathbf{a}}$ is good with probability $1-o(1)$.
\end{itemize}
\end{thm}

We now show why Theorem \ref{sieve-primes-2} implies Theorem \ref{sieve-primes}.  By \eqref{sigma-order}, we may choose $0 < c < 1/2$ small enough so that \eqref{sigma} holds. Let $A \geq 1$ be a fixed quantity.  Then we can find an integer $m$ obeying \eqref{moo} such that the quantity
$$ A' := 5^{-m} \times 80 c \log_2 x$$
is such that $A' \order A$ with implied constants independent of $A$.  

Suppose that we are in the probability $1-o(1)$ event that $\vec{\mathbf{a}}$ takes a value $\vec a$ which is good and such that \eqref{treat} holds.  On each sub-event $\vec{\mathbf{a}} = \vec a$ of this probability $1-o(1)$ event, we may apply Theorem \ref{packing-quant-cor} (for the random variables $\mathbf{n}_p$ conditioned to this event) define the random variables $\mathbf{n}'_p$ on this event with the stated properties.  For the remaining events $\vec{\mathbf{a}}=\vec a$, we set $\mathbf{n}'_p$ arbitrarily (e.g. we could set $\mathbf{n}'_p=0$).  The claim \eqref{up-short-random} then follows from Corollary \ref{packing-quant-cor} and \eqref{treat}, thus establishing Theorem \ref{sieve-primes}.

It remains to establish Theorem \ref{sieve-primes-2}.  This will be achieved in the next section.

\section{Using a sieve weight}\label{sec:weight}

If $r$ is a natural number, an \emph{admissible $r$-tuple} is a tuple $(h_1,\dots,h_r)$ of distinct integers $h_1,\dots,h_r$ that do not cover all residue
classes modulo $p$, for any prime $p$.  For instance, the tuple $(p_{\pi(r)+1},\dots,p_{\pi(r)+r})$ consisting of the first $r$ primes larger than $r$ is an admissible $r$-tuple.

We will establish Theorem \ref{sieve-primes-2} by a probabilistic argument involving a certain weight function.  More precisely, we will deduce this result from the following construction from \cite{FGKMT}.

\begin{thm}[Existence of good sieve weight]\label{weight}
Let $x$ be a sufficiently large real number, let $B_0$ be an integer, and let $y$ be any quantity obeying \eqref{ydef}.  Let $\PP, \QQ$ be defined by \eqref{p-def}, \eqref{q-def}.  Let $r$ be a positive integer with
\begin{equation}\label{r-bound}
r_0 \leq r \leq \log^{c_0} x
\end{equation}
for some sufficiently small absolute constant $c_0$ and sufficiently large absolute constant $r_0$, and let $(h_1,\dots,h_r)$ be an admissible $r$-tuple contained in $[2r^2]$.  Then one can find a positive quantity
\begin{equation}\label{alpha-crude}
\tau \geq x^{-o(1)}
\end{equation}
and a positive quantity $u = u(r)$ depending only on $r$ with
\begin{equation}\label{u-bound}
u \order \log r
\end{equation}
and a non-negative function $w: \PP \times \Z \to \R^+$ supported on $\PP \times (\Z \cap [-y,y])$ with the following properties:
\begin{itemize}
\item Uniformly for every $p \in \PP$, one has
\begin{equation}\label{wap}
 \sum_{n \in \Z} w(p,n) = \(1 + O\( \frac{1}{\log^{10}_2 x} \)\) \tau \frac{y}{\log^r x}.
\end{equation}
\item Uniformly for every $q \in \QQ$ and $i=1,\dots,r$, one has
\begin{equation}\label{wbp}
 \sum_{p \in \PP} w( p, q - h_i p ) = \(1 + O\( \frac{1}{\log^{10}_2 x} \)\) \tau \frac{u}{r} \frac{x}{2 \log^r x}.
\end{equation}
\item Uniformly for every $h = O(y/x)$ that is not equal to any of the $h_i$, one has
\begin{equation}\label{wcp}
\sum_{q \in \QQ} \sum_{p \in \PP} w( p, q - h p ) = O\( \frac{1}{\log^{10}_2 x} \tau \frac{x}{\log^r x} \frac{y}{\log x}\).
\end{equation}
\item Uniformly for all $p \in \PP$ and $n \in \Z$,
\begin{equation}\label{w-triv}
w(p,n) = O( x^{1/3+o(1)} ).
\end{equation}
\end{itemize}
\end{thm}

\begin{proof} See\footnote{The integer $B_0$ was not deleted from the sets $\PP$ or $\QQ$ in that theorem, however it is easy to see (using \eqref{w-triv}) that deleting at most one prime from either $\PP$ or $\QQ$ will not significantly worsen any of the estimates claimed by the theorem.} \cite[Theorem 5]{FGKMT}.  We remark that the construction of the weights and the verification of the required estimates relies heavily on the previous work of the second author in \cite{maynard-dense}.
\end{proof}

It remains to show how Theorem \ref{weight} implies Theorem \ref{sieve-primes-2}.  The analysis will be based on that in \cite[\S 5]{FGKMT}, which used a weight with slightly weaker hypotheses than in Theorem \ref{weight} to obtain somewhat weaker conclusions than Theorem \ref{sieve-primes-2} (in which the condition $q \equiv \mathbf{n}_p \pmod p$ was replaced by the stronger condition that $q = \mathbf{n}_p + h_i p$ for some $i=1,\dots,r$).

Let $x, B_0, c, y, z, \cS, \PP, \QQ$ be as in Theorem \ref{sieve-primes-2}.
Let $c_0$ be a sufficiently small absolute constant.  We set $r$ to be the maximum value permitted by Theorem \ref{weight}, namely
\begin{equation}\label{r-def}
 r := \lfloor \log^{c_0} x \rfloor 
\end{equation}
and let $(h_1,\dots,h_r)$ be the admissible $r$-tuple consisting of the first $r$ primes larger than $r$, thus $h_i = p_{\pi(r)+i}$ for $i=1,\dots,r$.  From the prime number theorem we have $h_i = O( r \log r )$ for $i=1,\dots,r$, and so we have $h_i \in [2r^2]$ for $i=1,\dots,r$ if $x$ is large enough (there are
many other choices possible, e.g. $(h_1,\ldots,h_r)=(1^2,3^2,\ldots,(2r-1)^2)$).  We now invoke Theorem \ref{weight} to obtain quantities $\tau,u$ and a weight $w: \PP \times \Z \to \R^+$ with the stated properties.

For each $p \in \PP$, let $\tilde{\mathbf{n}}_p$ denote the random integer with probability density
$$ \PR( \tilde{\mathbf{n}}_p = n ) := \frac{w(p,n)}{  \sum_{n' \in \Z} w(p,n') }$$
for all $n \in \Z$ (we will not need to impose any independence conditions on the $\tilde{\mathbf{n}}_p$).  From \eqref{wap}, \eqref{wbp} we have
\begin{equation}\label{wbp-diff}
\sum_{p \in \PP} \PR( q = \tilde{\mathbf{n}}_p + h_ip ) = \(1 + O\( \frac{1}{\log^{10}_2 x} \)\) \frac{u}{r} \frac{x}{2y}
\end{equation}
for every $q \in \QQ$ and $i=1,\dots,r$, and similarly from \eqref{wap}, \eqref{wcp} we have
\begin{equation}\label{wcp-diff}
\sum_{q \in \QQ} \sum_{p \in \PP} \PR( q = \tilde{\mathbf{n}}_p + hp ) \ll \frac{1}{\log^{10}_2 x} \frac{x}{\log x}
\end{equation}
for every $h = O(y/x)$ not equal to any of the $h_i$.  Finally, from \eqref{wap}, \eqref{w-triv}, \eqref{alpha-crude} one has
\begin{equation}\label{w-triv-diff}
 \PR( \tilde{\mathbf{n}}_p = n ) \ll x^{-1/2 - 1/6 +o(1)}
\end{equation}
for all $p \in \PP$ and $n \in \Z$.

We choose the random vector
$\vec{\mathbf{a}} := (\mathbf{a}_s \mod s)_{s \in \cS}$ by selecting each $\mathbf{a}_s \mod s$ uniformly at random from $\Z/s\Z$, independently in $s$ and independently of the $\tilde{\mathbf{n}}_p$.
The resulting sifted set $S(\vec{\mathbf{a}})$ is a random periodic subset of $\Z$ with density
$$ \sigma := \prod_{s \in \cS} \(1 - \frac{1}{s}\).$$
From the prime number theorem (with sufficiently strong error term), 
\eqref{zdef} and \eqref{s-def},
\[
\sigma = \(1 + O\(\frac{1}{\log^{10}_2 x}\)\) \frac{\log(\log^{20} x)}{\log z} \\
= \(1 + O\(\frac{1}{\log^{10}_2 x}\)\) \frac{80 \log_2 x}{\log x \log_3 x / \log_2 x},
\]
so in particular we see from \eqref{ydef} that
\begin{equation}\label{gamma-y}
\sigma y = \(1 + O\(\frac{1}{\log^{10}_2 x}\)\) 80 c x \log_2 x.
\end{equation}
We also see from \eqref{r-def} that
\begin{equation}\label{gamma-small}
\sigma^r = x^{o(1)}.
\end{equation}

We have a useful correlation bound:

\begin{lem}\label{gamma-cor}  Let $t \le \log x$ be a natural number, and let $n_1,\dots,n_t$ be distinct integers of magnitude $O(x^{O(1)})$.  Then one has
$$ \PR( n_1,\dots,n_t \in S(\vec{\mathbf{a}}) ) = \(1 + O\pfrac{1}{\log^{16} x}\) \sigma^t.$$
\end{lem}

\begin{proof}  See \cite[Lemma 5.1]{FGKMT}.
\end{proof}

Among other things, this gives the claim \eqref{treat}:
 
\begin{cor}\label{s0}  For any fixed $0 \leq \alpha < \beta \leq 1$, we have with probability $1-o(1)$ that
\begin{equation}\label{qqa}
\# (\QQ \cap [\alpha y, \beta y] \cap S(\vec{\mathbf{a}}) ) \sim \sigma |\beta-\alpha| \frac{y}{\log x} \sim
 80 c |\beta-\alpha| \frac{x}{\log x} \log_2 x.
\end{equation}
\end{cor}

\begin{proof}  See \cite[Corollary 4]{FGKMT}, replacing $\QQ$ with $\QQ \cap [\alpha y, \beta y]$.
\end{proof}

For each $p \in \PP$, we consider the quantity
\begin{equation}\label{xp-def}
 X_p(\vec{a}) := \PR( \tilde{\mathbf{n}}_p + h_ip \in S(\vec{a}) \text{ for all } i=1,\dots,r ),
\end{equation}
and let $\PP(\vec a)$ denote the set of all the primes $p \in \PP$ such that
\begin{equation}\label{sumn}
X_p(\vec a) = \(1 + O_{\leq}\(\frac{1}{\log^3 x}\)\) \sigma^r.
\end{equation}

In light of Lemma \ref{gamma-cor}, we expect most primes
in $\PP$ to lie in $\PP(\vec{a})$, and this will be confirmed
below (Lemma \ref{smc}).
We now define the random variables $\mathbf{n}_p$ as follows.  Suppose we are in the event $\vec{\mathbf{a}} = \vec a$ for some $\vec a$ in the range of $\vec{\mathbf{a}}$.  If $p \in \PP \backslash \PP(\vec{a})$, we set $\mathbf{n}_p=0$.  Otherwise, if $p \in \PP(\vec{a})$, we define $\mathbf{n}_p$ to be the random integer with conditional probability distribution
\begin{equation}\label{xpa}
 \PR( \mathbf{n}_p = n | \vec{\mathbf{a}} = \vec a ) := \frac{Z_p(\vec{a};n)}
{X_p(\vec a)}, \quad Z_p(\vec{a};n) = 1_{n+h_jp\in
S(\vec{a})\text{ for }j=1,\ldots,r} \PR(\tilde{\mathbf{n}}_p=n).
\end{equation}
with the $\mathbf{n}_p$ jointly conditionally independent on the event $\vec{\mathbf{a}} = \vec a$.  From \eqref{xp-def} we see that these random variables are well defined.

Substituting definition \eqref{xpa} into the left hand side of \eqref{good},
and observing that $\mathbf{n}_p\equiv q\pmod{p}$ is only possible if 
$p\in \PP(\vec{\ba})$, we see that to prove \eqref{good}, it suffices to show
 that with probability 
$1-o(1)$ in $\vec{\ba}$, for all but at most $\frac{x}{\log x\log_2 x}$ primes in
$\QQ \cap S(\vec{\ba})$, we have
\be\label{good-1}
\sigma^{-r}
\sum_{p\in\PP(\vec{\ba})} \sum_{h} Z_p(\vec{\ba};q-hp) = C +O\pfrac{1}{\log_2^3 x}.
\ee

We now confirm that $\PP \backslash \PP(\vec{\ba})$ is small with high
probability.

\begin{lem}\label{smc}   
With probability $1-O(1/\log^3 x)$, $\PP(\vec{\mathbf{a}})$ contains
 all but $O( \frac{1}{\log^3 x} \frac{x}{\log x})$ of the primes $p \in \PP$. 
In particular, $\E \# \PP(\vec{\mathbf{a}}) = \# \PP (1+O(1/\log^3 x))$.
\end{lem}

\begin{proof}  See \cite[Lemma 5.3]{FGKMT}.
\end{proof}

The left side of relation \eqref{good-1} breaks naturally into two pieces,
a `main term' consisting of summands where $h=h_i$ for some $i$, and 
an `error terms' consisting of the remaining summands.  We first take care of 
the error terms.

\begin{lem}\label{smc-1}
With probability $1-o(1)$ we have
\be\label{good-h}
\sigma^{-r} \sum_{p\in \PP(\vec{\ba})} \sum_{\substack{h\ll y/x\\ h\not\in\{h_1,\ldots,h_r\}}}  
Z_p(\vec{\ba};q-hp) \ll \frac{1}{\log_2^3 x}
\ee
for all but at most $\frac{x}{2\log x\log_2 x}$ primes 
$q\in \QQ\cap S(\vec{\ba})$.
\end{lem}

\begin{proof}
We first extend the sum over all $p\in\PP$. By Markov's inequality, it suffices to show that
\be\label{ep} \E \  \sum_{q \in \QQ \cap S(\vec{\mathbf{a}}) } \sigma^{-r} \sum_{p \in \PP}\sum_{\substack{h\ll y/x\\ h\notin \{h_1,\dots,h_k\}}}
Z_p(\vec{\ba};q-hp)
=o\( \frac{x}{\log x\log^{4}_2 x}\).\ee
The left-hand side of \eqref{ep} equals
$$
\sigma^{-r}\sum_{q \in \QQ}  \sum_{\substack{h\ll y/x\\ h\notin \{h_1,\dots,h_k\}}}  \sum_{p \in \PP}
\PR(q\in S(\vec{\ba}), q+h_jp-hp\in S(\vec{\ba})\text{ for }j=1,\ldots,r)
\PR(q=\tilde{\mathbf{n}}_p+hp).$$
We note that for any $h$ in the above sum, the $r+1$ integers $q, q+h_1 p-hp, \dots, q + h_r p-hp$ are distinct.
Applying Lemma \ref{gamma-cor}, followed by \eqref{wcp-diff}, we may thus bound
this expression by
$$ \ll  \sum_{\substack{h\ll y/x\\ h\notin\{h_1,\dots,h_k\}}} \sigma \ \frac{x/\log x}{\log_2^{10} x} \ll
 \sigma \frac{1}{\log^{10}_2 x} \frac{y}{\log x}.$$
The claim now follows from \eqref{gamma-y}.
\end{proof}

Next, we deal with the main term of \eqref{good-1},
by showing an analogue of \eqref{wbp-diff}.

\begin{lem}\label{smc-2}  
With probability $1-o(1)$, we have
\begin{equation}\label{sumno}
\sigma^{-r} \sum_{i=1}^r \sum_{p\in \PP(\vec{\ba})} 
Z_p(\vec{\ba};q-h_ip) = \(1 + O\(\frac{1}{\log^{3}_2 x}\)\) 
\frac{u}{\sigma} \frac{x}{2y}
\end{equation}
for all but at most $\frac{x}{2\log x \log_2 x}$ of the primes $q \in \QQ \cap S(\vec{\mathbf{a}})$.
\end{lem}

\begin{proof} 
We first show that replacing $\PP(\vec{\ba})$ with $\PP$
has negligible effect on the sum, with probability $1-o(1)$.  
Fix $i$ and susbtitute $n=q-h_ip$.  By Markov's inequality, it
suffices to show that
\be\label{soo-2}
\E \sum_n \sigma^{-r} \sum_{p\in \PP\backslash \PP(\vec{\ba})} Z_p(\vec{\ba};n)
=o\(
\frac{u}{\sigma} \frac{x}{2y}\ \frac{1}{r} \frac{1}{\log_2^3 x}\ \frac{x}{\log x\log_2 x} \).
\ee
By Lemma \ref{gamma-cor}, we have
\begin{align*}
\E\ \sum_{n} \sigma^{-r} \sum_{p \in \PP} Z_p(\vec{\ba};n)
&= \sigma^{-r} \sum_{p \in \PP} \sum_{n}  \PR(\tilde{\mathbf{n}}_p=n)
\PR( n+h_jp\in S(\vec{\ba})\text{ for }j=1,\ldots,r ) \\
&=\(1+O\pfrac{1}{\log^{16} x}\) \# \PP.
\end{align*}
Next, by \eqref{sumn} and Lemma \ref{smc} we have
\begin{align*}
\E\ \sum_{n} \sigma^{-r} &\sum_{p \in \PP(\vec{\ba})}  Z_p(\vec{\ba};n)
= \sigma^{-r} \sum_{\vec{a}} \PR(\vec{\ba}=\vec{a}) \sum_{p\in \PP(\vec{a})} 
X_p(\vec{a}) \\
&=\(1 + O\pfrac{1}{\log^3 x}\) \ 
\E \; \# \PP(\vec{\mathbf{a}}) 
=\(1 + O\pfrac{1}{\log^3 x}\) \# \PP ;
\end{align*}
subtracting, we conclude that the left-hand side of \eqref{soo-2} is
$O(\#\PP/\log^3 x)=O(x/\log^4 x)$. 
The claim then follows from \eqref{ydef} and \eqref{r-bound}.

By \eqref{soo-2}, it suffices to show that with probability $1-o(1)$,
for all but at most $\frac{x}{2\log x\log_2 x}$ primes $q \in \QQ \cap S(\vec{\mathbf{a}})$, one has
\begin{equation}\label{soo}
\sum_{i=1}^r \sum_{p\in \PP} Z_p(\vec{\ba};q-h_ip)
 = \(1 + O_{\le}\(\frac{1}{\log^{3}_2 x}\)\) 
\sigma^{r-1} u \frac{x}{2y}.
\end{equation}

Call a prime $q \in \QQ$ \emph{bad} if $q \in \QQ \cap S(\vec{\mathbf{a}})$ but
\eqref{soo} fails.  
Using Lemma \ref{gamma-cor} and \eqref{wbp-diff}, we have
\begin{align*}
\E \bigg[\sum_{q\in\QQ\cap S(\vec{\mathbf{a}}) } \sum_{i=1}^r \sum_{p \in \PP} Z_p(\vec{\ba};q-h_ip)
\bigg]&=\sum_{q,i,p} \PR( q + (h_j-h_i) p \in S(\vec{\mathbf{a}}) \text{ for all } j=1,\dots,r)\PR(\tilde{\mathbf{n}}_p=q-h_ip)\\
&=\(1+O\pfrac{1}{\log_2^{10} x}\) \frac{\sigma y}{\log x} \ 
\sigma^{r-1} u \ \frac{x}{2y}
\end{align*}
and
\begin{align*}
\E \bigg[ \sum_{q\in\QQ\cap S(\vec{\mathbf{a}}) } \bigg( \sum_{i=1}^r \sum_{p\in\PP} Z_p (\vec{\ba};q-h_ip)\bigg)^2
\bigg]&=\sum_{\substack{p_1,p_2,q \\ i_1,i_2}} \PR( q + (h_j-h_{i_\ell}) p_\ell \in S(\vec{\mathbf{a}}) 
\text{ for } j=1,\dots,r;\ell=1,2) \\
& \qquad \times \PR(\tilde{\mathbf{n}}^{(1)}_{p_1}=q-h_{i_1}p_1) 
\PR(\tilde{\mathbf{n}}^{(2)}_{p_2}=q-h_{i_2}p_2) \\
&=\(1+O\pfrac{1}{\log_2^{10} x}\)  \frac{\sigma y}{\log x} \ 
\(\sigma^{r-1} u \ \frac{x}{2y}\)^2,
\end{align*}
where $(\tilde{\mathbf{n}}^{(1)}_{p_1})_{p_1 \in \PP}$ and $(\tilde{\mathbf{n}}^{(2)}_{p_2})_{p_2 \in \PP}$ are independent copies of $(\tilde{\mathbf{n}}_p)_{p \in \PP}$ over $\vec{\mathbf{a}}$.  In the last step we used the fact that the terms with
$p_1=p_2$ contribute negligibly.

By Chebyshev's inequality (Lemma \ref{cheb}) it follows that the number of
bad $q$ is $\ll \frac{\sigma y}{\log x} \frac{1}{\log_2^3 x} \ll \frac{x}{\log x\log_2^2 x}$ with probability $1-O(1/\log_2 x)$.  This concludes the proof.
\end{proof}

We now conclude the proof of
Theorem \ref{sieve-primes-2}.  We need to prove \eqref{good-1};
this follows immediately from Lemma \ref{smc-1} and Lemma \ref{smc-2}
 upon noting that by \eqref{r-def}, \eqref{u-bound} and \eqref{gamma-y},
$$
C := \frac{u}{\sigma} \ \frac{x}{2y} \asym \frac{1}{c}.
$$

\end{document}